\documentclass[microtype]{gtpart}
\usepackage{pinlabel}

\title{Stable systolic category of the product of spheres}

\author{Hoil Ryu}
\givenname{Hoil}
\surname{Ryu}
\address{Graduate School of Mathematics\\Kyushu University\\\newline
	774\\Motooka\\Nishi-ku\\Fukuoka\\819-0395\\Japan}
\email{h-ryu@math.kyushu-u.ac.jp}
\urladdr{}


\arxivreference{}
\arxivpassword{}

\volumenumber{}
\issuenumber{}
\publicationyear{}
\papernumber{}
\startpage{}
\endpage{}
\doi{}
\MR{}
\Zbl{}
\received{}
\revised{}
\accepted{}
\published{}
\publishedonline{}
\proposed{}
\seconded{}
\corresponding{}
\editor{}
\version{}

\newtheorem{thm}{Theorem}[section]
\newtheorem{cor}{Corollary}[section]
\newtheorem{lemma}{Lemma}[section]

\newtheorem{prop}{Proposition}[section]

\newtheorem*{GromovThm}{Gromov's Theorem}

\theoremstyle{definition}
\newtheorem{defn}{Definition}[section]
\newtheorem*{rmk}{Remark}
\newtheorem*{example}{Example}

\makeatletter
\let\c@cor=\c@thm
\let\c@lemma=\c@thm
\let\c@claim=\c@thm
\let\c@prop=\c@thm
\let\c@conjecture=\c@thm
\let\c@prob=\c@thm
\let\c@defn=\c@thm
\makeatother

\makeautorefname{thm}{Theorem}
\makeautorefname{cor}{Corollary}
\makeautorefname{lemma}{Lemma}
\makeautorefname{claim}{Claim}
\makeautorefname{prop}{Proposition}
\makeautorefname{conjecture}{Conjecture}
\makeautorefname{prob}{Problem}

\makeautorefname{GromovThm}{Gromov's Theorem}

\newcounter{fig}

\DeclareMathOperator{\support}{supp}

\DeclareMathOperator{\interior}{Int}

\DeclareMathOperator{\Lip}{Lip}

\DeclareMathOperator{\mass}{mass}
\DeclareMathOperator{\stmass}{stmass}
\DeclareMathOperator{\comass}{comass}

\newcommand{\current}{\mathcal{D}}

\newcommand{\polyhedralD}{\mathcal{P}}
\newcommand{\normalD}{\mathcal{N}}
\newcommand{\rectifiableD}{\mathcal{R}}
\newcommand{\flatD}{\mathcal{F}}

\DeclareMathOperator{\stsys}{stsys}
\DeclareMathOperator{\catstsys}{cat_{stsys}}
\DeclareMathOperator{\cuplength}{cup}

\DeclareMathOperator{\size}{size}

\DeclareMathOperator{\LPD}{lpd}
\DeclareMathOperator{\MOD}{mod}

\begin{document}

\begin{abstract}
The stable systolic category of a closed manifold $M$ indicates the complexity in the sense of volume.
This is a homotopy invariant, even though it is defined by some relations between homological volumes on $M$.
We show an equality of the stable systolic category and the real cup-length for the product of arbitrary finite dimensional real homology spheres.
Also we prove the invariance of the stable systolic category under the rational equivalences for orientable $0$--universal manifolds.
\end{abstract}

\begin{asciiabstract}
The stable systolic category of a closed manifold M indicates the complexity in the sense of volume.
This is a homotopy invariant, even though it is defined by some relations between homological volumes on M.
We show an equality of the stable systolic category and the real cup-length for the product of arbitrary finite dimensional real homology spheres.
Also we prove the invariance of the stable systolic category under the rational equivalences for orientable 0-universal manifolds.
\end{asciiabstract}

\maketitle


\section{Introduction}

In this paper, a manifold is assumed to be closed, connected, orientable and smooth.
The systole of a manifold $M$ is the least length of non-contractible closed loops in $M$.
One can generalize this concept to the least volume of $k$--dimensional nonzero homology classes, so called as the homology systole.
Now we can imagine such systoles have some kind of relations with the entire volume of $M$, 
and it is natural to ask what kind of relationship exists.

As an answer, Gromov proved a theorem which says that the existence of non-trivial cup product implies the existence of the stable isosystolic inequality as follows.

\begin{GromovThm}{\rm \cite[7.4.C]{Gro83}}\qua
\label{thm:Gromov83}
Let $M$ be an $n$--manifold (which can be non-orientable).
If there exist some reduced real cohomology classes $\alpha_1^*, \cdots, \alpha_k^*$ with $\alpha_i^*$ in $\tilde{H}^{d_i}(M;\mathbb{R})$ and a nonzero cup product $\alpha_1^* \cup \cdots \cup \alpha_k^*$ in $\tilde{H}^n(M;\mathbb{R})$,
then there exists $C > 0$ satisfying
\begin{align*}
	\prod_{i=1}^{k}\stsys_{d_i}(M,\mathcal{G}) \le C \cdot \mass\bigl( [M], \mathcal{G} \bigr)
\end{align*}
for all Riemannian metric $\mathcal{G}$ on $M$
where $\stsys_{d_i}$ is the stable $d_i$--systole and $[M]$ is the fundamental class of $M$ with coefficients in $\mathbb{Z}$ for orientable cases or $\mathbb{Z}/2\mathbb{Z}$ for non-orientable cases.
\end{GromovThm}

The greatest $k$ satisfying the stable isosystolic inequality is called as the stable systolic category of $M$ which is introduced by Katz and Rudyak \cite{KatRud06}, and it is known as a homotopy invariant by Katz and Rudyak \cite{KatRud08}.
We will show the stable systolic category of $0$-universal manifold is also invariant under the rational equivalences in \fullref{cor:REinvariance}.

For an orientable manifold $M$, Gromov's Theorem implies that the stable systolic category is not smaller than the real cup-length.
So, is there some manifold $M$ such that the stable systolic category is greater than the real cup-length?
If such $M$ exists, then the inversion of Gromov's Theorem will fail for $M$, while
this interesting question is not answered yet.
Instead of the answer, it is known the equality of them for some manifolds, eg,
Dranishnikov and Rudyak \cite{DraRud09}.
In this paper, we also show more equality later in
\fullref{thm:LSSCPstsyscat} and \fullref{thm:CatstsysProductSphere}.

\subsection{Definition of the stable systolic category}

Let $(X,A)$ be an object of the local Lipschitz category $\mathfrak{L}$.
In this paper, we suppose $X$ is a nonempty subset of some finite dimensional Euclidean space $\mathbb{R}^n$ with the standard norm.
So $X$ possesses the restricted metric of $\mathbb{R}^n$.
Let $G$ be a $\mathbb{Z}$--module with a norm $| \cdot |$ which makes $G$ a complete metric space.
If $G$ is $\mathbb{Z}$ or $\mathbb{R}$, we assume that norm of $G$ is the standard norm.
The comass of a differential form $\omega$ on $X$ is defined as
$$
	\comass(\omega) := \sup \bigl\{ | \omega_x(\tau) | : x \in X, \text{orthonormal } q \text{--frame } \tau \bigr\} .
$$
Also, the mass of a $q$--current $T$ in $X$ is the dual norm of comass, ie,
$$
	\mass(T) := \sup \bigl\{ T(\omega) : \text{differential } q\text{--form } \omega, \comass(\omega) \le 1 \bigr\}.
$$
A Lipschitzian singular $q$--cube $\kappa\co I^q \to X$,
induces a homomorphism $\kappa_{\flat}$ from the module of polyhedral chains $\polyhedralD_q(X;G)$ to the module of rectifiable currents $\rectifiableD_q(X;G)$.
Then the mass of $\kappa$ is defined by the mass of the image $\kappa_{\flat}I^q$
where $I^q$ is the corresponding polyhedral $q$--current of the unit rectangular parallelepiped $I^q$.
This correspondence of $\kappa$ to $\kappa_{\flat}I^q$ gives a chain map $\Phi$ of degree $0$ from the chain complex of all Lipschitzian singular cubes into the chain complex of flat chains $\flatD_*(E|X)$. 
Then we can verify that $\Phi$ induces an isomorphism $\Phi_*$ between the flat homology module $H^{\flat}_q(X,A;G)$ and the singular homology module $H_q(X,A;G)$ for all $q$.
One can find more details of these definitions at Federer \cite{Fed69}, Federer \cite{Fed74}, Federer and Fleming \cite{FedFle60}, Serre \cite{Ser51} and White \cite{White99}.

For a Lipschitzian singular chain $c$, there exists a representation $\sum_i \kappa_i \otimes g_i$ where $g_i$ is contained in $G$ and $\kappa_i$ is a Lipschitzian singular $q$--cube which is not overlapping each other (subdivide if necessary).
Then the {\it mass} of $c$ is defined as
$$
	\mass(c) := \sum_i |g_i| \cdot \mass(\kappa_i) \,.
$$
%
The {\it mass} or {\it volume} of a singular homology class $\eta$ in $H_q(X,A;G)$ is defined by
$$
	\mass(\eta;G) := \inf\bigl\{ \mass(c) : \eta = [c],\ c \text{ is a Lipschitzian cycle} \bigr\} \,.
$$
If $G$ is $\mathbb{R}$, the mass is a norm on the homology vector spaces.
We will omit $G$ in the case of $\mathbb{Z}$.


The {\it $q$--dimensional homology systole} of $(X,A)$ is defined by infimum of mass of non-trivial $q$--th integral homology classes.
However Gromov \cite{Gro83} claims that Gromov's Theorem will fail for $S^1 \times S^3$, if we consider the homology systoles instead of the stable systoles.
Briefly, we can consider the stable systole as a systole in the real homology vector spaces.
Here we give formal definition for the stable systole.
The {\it stable mass} on $H_q(X,A;\mathbb{Z})$ is defined as 
$$
	\stmass(\eta) := \inf\bigl\{ \mass(m \cdot \eta) / m : 0 < m \in \mathbb{Z} \bigr\}
$$
for all $\eta \in H_q(X,A;\mathbb{Z})$.
The inclusion $\iota\co \mathbb{Z} \to \mathbb{R}$ induces the coefficient homomorphism $\iota_*$ on homology.
Federer \cite[5.8]{Fed74} showed the stable mass of $\eta$ is equal to the mass of the image $\iota_*\eta$.
So we can define the {\it $q$--dimensional stable systole} of $(X,A)$ as
\begin{align*}
	\stsys_q(X,A) :=
	\inf\left\{ \stmass(\eta) : \eta \in H_q(X,A;\mathbb{Z}),\ \iota_*\eta \neq 0 \right\}.
\end{align*}
A homology $q$--systole or a stable $q$--systole is called {\it trivial}, if it is infinite.
If the $q$--th real homology vector space $H_q(X,A;\mathbb{R})$ is zero, then the stable $q$--systole is trivial for all Riemannian metrics on $(X,A)$.
Hence if the $q$--th integral homology module $H_q(X,A;\mathbb{Z})$ is a torsion module, then the stable $q$--systole is trivial for every metric on $(X,A)$.


For a given positive integer $n > 0$, 
a $k$--tuple $P = (p_1, \cdots, p_k)$ of positive integers is called 
a {\it partition} of $n$ if
$n = p_1 + \cdots + p_k$ and $p_1 \le \cdots \le p_k \le n$.
A partition $P$ is called {\it positive} (or {\it non-negative}) if $p_i > 0$ (or $p_i \ge 0$) for all $i$.
The {\it size} of a partition which denoted by $\size(P)$ is defined by the cardinality of positive integers contained in the partition.
Hence if a $k$--tuple $P$ is a positive partition, then the size of partition is $k$.
From now on, we suppose a partition is positive unless otherwise stated.
For a partition $P$,
the {\it duplicated number} of $p_i$ is the cardinality of elements in $P$ who are equal to $p_i$.

Now we define concepts for an $n$--manifold $M$.
A partition $P$ of $n$ is called {\it stable systolic categorical} for $M$, if there exists a real number $C > 0$ and non-trivial stable $p_i$--systoles such that
$$\prod_{i=1}^{\size(P)}\stsys_{p_i}(M,\mathcal{G}) \le C \cdot \mass\bigl( [M], \mathcal{G}; \mathbb{Z}/2\mathbb{Z} \bigr)$$
for every Riemannian metric $\mathcal{G}$ on $M$ where the fundamental class $[M]$ in $H_n(M;\mathbb{Z}/2\mathbb{Z})$.

\begin{defn}
The {\it stable systolic category} of $M$ is defined by
$$\catstsys(M):=\sup \bigl\{\mathrm{size}(P):P \text{ is stable systolic categorical partition for } M \bigr\}.$$
\end{defn}

The {\it real cup-length} of $M$ is defined by
$$\cuplength_{\mathbb{R}}(M) := \min \big\{ k \ge 0 : \alpha_0 \cup \alpha_1 \cup \cdots \cup \alpha_k = 0 \text{ for all } \alpha_i \in \tilde{H}^*(M;\mathbb{R}) \big\}$$
where $\tilde{H}^*(M;\mathbb{R})$ is the reduced real cohomology ring of $M$.
As we said before, the real cup-length is a lower estimate for the stable systolic category from Gromov's Theorem.

If $M$ is non-orientable, then the top dimensional real cohomology vector space $H^n(M;\mathbb{R})$ vanishes.
So every cohomology class in $H^n(M;\mathbb{R})$ vanishes, we can not apply Gromov's Theorem for top dimension.
This is a reason to consider only orientable manifolds in this paper.

\subsection{Acknowledgments}

The author expresses gratitude to Professor Norio Iwase, whose guidance and support to write this paper.

\section{Preliminaries on the stable systoles}


Many equations and inequalities for mass are studied.
One can find those results at Babenko \cite{Bab93}, Federer \cite{Fed69} and Whitney \cite{Whi57}.
Here we state or recall some of them for the stable systoles, with some appropriate modifications applied.

\begin{prop}
\label{prop:stable0systole}
For a local Lipschitz neighborhood retract $X$ in $\mathbb{R}^n$,
the stable $0$--systole is $1$.
\end{prop}

\begin{proof}
Let $\current_0(X)$ be the vector space of $0$--currents.
A map $\mathfrak{d}\co X \to \current_0(X)$ can be defined as
$\mathfrak{d}(x)(\omega) = \mathfrak{d}_{x}(\omega) := \omega(x)$ for a point $x$ of $X$ and a differential $0$--form $\omega$ on $X$.
Then $\mathfrak{d}_x$ is a polyhedral $0$--current with $\mass(\mathfrak{d}_x) = 1$.
This implies that $\mathfrak{d}_{x}$ is a normal $0$--cycle with coefficients~$\mathbb{Z}$.
Furthermore, the image $\iota_*\Phi_*^{-1}[\mathfrak{d}_x]$ is not vanished in $H_0(X;\mathbb{R})$.
So we have
$$\stsys_0(X) = \mass\bigl( \iota_*\Phi_*^{-1}[\mathfrak{d}_x] \bigr) = 1$$
for an arbitrary point $x$ in $X$.
\end{proof}

\begin{lemma}
\label{lemma:MassRescaleEquality}
For a local Lipschitz neighborhood retract $X$ in $\mathbb{R}^n$,
if one rescale the standard metric $\mathcal{G}$ on $\mathbb{R}^n$ by the square of a real number $t > 0$, 
then the quotient mass of a homology class $\eta \in H_q(X;G)$ increase by the $t^q$ times.
Furthermore, the stable $q$--systole satisfies
$$\stsys_q(X,t^2 \mathcal{G}|X) = t^q \cdot \stsys_q(X, \mathcal{G}|X)$$
where $\mathcal{G}|X$ is the restriction of $\mathcal{G}$ on $X$.
\end{lemma}

\begin{proof}
A similar result was introduced by Whitney \cite{Whi57} for the real flat chains.
So the first result is satisfied for an arbitrary homology class.
Also the definition of the stable systole implies
\begin{align*}
	\stsys_q(X,t^2 \mathcal{G}|X)
	& = \inf\left\{ t^q \cdot \mass(\iota_*\eta, \mathcal{G}|X; \mathbb{R}) : \eta \in H_q(X,A;\mathbb{Z}),\ \iota_*\eta \neq 0 \right\} 
\end{align*}
which means the equality for the stable systoles.
\end{proof}


\begin{prop}{\cite[X.6 and X.7]{Whi57}}\qua
\label{prop:PushForwardMassInequality}
Let $X$ and $Y$ be open subsets in $\mathbb{R}^m$ and $\mathbb{R}^n$ respectively.
For a locally Lipschitzian map $f\co X \to Y$ and an integral rectifiable $q$--current $T$ whose support is contained in a compact subset $K$ of $X$,
there exists an inequality 
$$\mass(f_{\flat}T) \le \Lip(f|K)^q \cdot \mass(T)$$
where $\Lip(f|K)$ is the lower bound of Lipschitz constants of the restriction $f|K$.
\end{prop}

\begin{prop}
\label{prop:MappedLowerBoundMass}
Let $(X,A)$ and $(Y,B)$ be local Lipschitz neighborhood retract pairs in $\mathbb{R}^m$ and $\mathbb{R}^n$ respectively.
If $f\co (X,A) \to (Y,B)$ is a locally Lipschitzian map,
then for any homology class $\eta$ of $H_q(X,A;G)\,$,
there is a compact subset $K$ of $\mathbb{R}^m$ which satisfies
$$0 \le \mass(f_*\eta;G) \le \Lip(f|K)^q \cdot \mass(\eta;G)$$
where $f_*\co H_q(X,A;G) \to H_q(Y,B;G)$ is the induced homomorphism.
\end{prop}

\begin{proof}
Note that $f$ induces a homomorphism $f_\flat\co Z_q(X,A;G) \to Z_q(Y,B;G)$ on flat cycles as well as $f_{\flat}\flatD_q(\mathbb{R}^m|A;G) \subset \flatD_q(\mathbb{R}^n|B;G)$ .
For a given flat homology class $\Phi_*\eta$,
let $T$ be a representative normal $q$--cycle in $Z_q(X,A;G)$.
The naturality of $\Phi_*$ implies $\Phi_* f_* \eta = f_* \Phi_* \eta = f_*[T] = [f_{\flat} T]$.
Also the relation of cosets
$[f_{\flat} T] = [f_{\flat}T+f_{\flat}\flatD_q(\mathbb{R}^m|A;G)] = [f_{\flat}T+\flatD_q(\mathbb{R}^n|B;G)]$
implies that the relation of the sets
$$\bigl\{ f_\flat T : [T] = \Phi_*\eta \bigr\} \subset \bigl\{ S : [S] = \Phi_*f_*\eta \bigr\} \subset Z_q(Y,B;G) \,.$$
With the definition of the mass of homology class, we obtain
\begin{align*}
	\mass(f_*\eta ; G)
	& \le \inf\bigl\{ \mass(f_\flat T) : [T] = \Phi_*\eta \bigr\}.
\end{align*}
Because of $T$ is compact supported, there is a compact subset $K$ of $\mathbb{R}^m$ with $\support(T) \subset \interior(K)$.
Here we can apply \fullref{prop:PushForwardMassInequality} for $T$, so we have
\begin{align*}
	\mass(f_*\eta ; G)
	& \le \Lip(f|K)^q \cdot \inf\bigl\{ \mass(T) : [T] = \Phi_*\eta \bigr\} 
\end{align*}
which implies the result.
\end{proof}

\begin{lemma}
\label{lemma:StableSystoleMappedInequality}
Let $(X,A)$ and $(Y,B)$ are local Lipschitz neighborhood retract pairs.
If a locally Lipschitzian map $f\co (X,A) \to (Y,B)$ induces
a monomorphism $f_*\co H_q(X,A;\mathbb{R}) \to H_q(Y,B;\mathbb{R})$,
then there is a compact subset $K$ in the ambient space of $X$ satisfying
$$\stsys_q(Y,B) \le \Lip(f|K)^q \cdot \stsys_q(X,A) .$$
Furthermore, if $H_q(X,A;\mathbb{R})$ is nonzero, then $\stsys_q(Y,B)$ is a positive real number.
\end{lemma}

\begin{proof}
\fullref{prop:MappedLowerBoundMass} and $f_*\bigl(H_q(X,A;\mathbb{R}) \setminus \{0\}\bigr) \subset \bigl(H_q(Y,B;\mathbb{R}) \setminus \{0\}\bigr)$
imply the existence of inequality in the stable systole level.

For integral homology class $\eta$ with $\iota_*\eta$ is nonzero,
the image $f_*\iota_*\eta$ does not vanish, since $f_*$ is a monomorphism.
Recall that the mass of real homology classes is a norm, hence $\mass(f_*\iota_*\eta)$ is a positive real number.
Furthermore, the stable $q$--systole does not converges to zero, since $\mathbb{Z}$ is discrete.
\end{proof}


\begin{prop}
\label{prop:MassCrossEquality}
Let $X$ and $Y$ be open subsets of $\mathbb{R}^m$ and $\mathbb{R}^n$ respectively.
For rectifiable currents $S$ in $\rectifiableD_p(X)$ and $T$ in $\rectifiableD_q(Y)$,
the mass of their cross product is equal to the multiplication of their masses, ie,
$$\mass(S \times T) = \mass(S) \cdot \mass(T)$$
with respect to the product metric on $X \times Y$.
\end{prop}

\begin{proof}
Since $S$ and $T$ are rectifiable currents, mass can be written by associated Radon measures $\Vert S \Vert$, $\Vert T \Vert$ and $\Vert S \times T \Vert$.
Therefore Fubini's Theorem (see Federer \cite[2.6.2.(2)]{Fed69}) implies
$$\mass(S \times T) 
= \Vert S \times T \Vert(X \times Y) 
= \Vert S \Vert (X) \cdot \Vert T \Vert (Y)  
= \mass(S) \cdot \mass(T)$$
the result.
\end{proof}

\begin{lemma}
\label{lemma:QuotientMassCrossEquality}
Let $(X,A)$ and $(Y,B)$ are local Lipschitz neighborhood retract pairs.
For homology classes $\xi \in H_p(X,A;G)$ and $\eta \in H_q(Y,B;G)$,
we can estimate
\begin{align*}
	\mass(\xi \times \eta; G) & \le \mass(\xi; G) \cdot \mass(\eta; G) \\
\tag*{\text{\sl and}}
	\stsys_{p+q}\bigl( (X,A) \times (Y,B) \bigr) & \le \stsys_p(X,A) \cdot \stsys_q(Y,B)
\end{align*}
with respect to the product metric on $(X,A) \times (Y,B)$.
\end{lemma}

\begin{proof}
Let $S$ and $T$ be representative rectifiable cycles corresponding to $\xi$ and $\eta$ respectively, ie, $\Phi_*\xi = [S]$ with $S \in Z^{\flat}_{p}(X,A;G)$ and $\Phi_*\eta = [T]$ with $T \in Z^{\flat}_{q}(Y,B;G)$.
Then the naturality of a cross product implies that there is a representative rectifiable current with the form of a cross product $S \times T$ in the coset $[c] = \Phi_*(\xi \times \eta)$. Therefore
\begin{align*}
	\bigl\{ S \times T : [S] \times [T] = \Phi_*\xi \times \Phi_*\eta \bigr\}
	& = \bigl\{ S \times T : [S \times T] = \Phi_*(\xi \times \eta) \bigr\} \\
	& \subset \bigl\{ c : [c] = \Phi_*(\xi \times \eta) \bigr\} \\
	& \subset Z^{\flat}_{p+q}\bigl( (X,A) \times (Y,B) ; G \bigr) .
\end{align*}
Hence \fullref{prop:MassCrossEquality} implies an inequality
\begin{align*}
	\mass(\xi \times \eta; G)
	& \le \inf\bigl\{ \mass(S \times T) : [S] \times [T] = \Phi_*\xi \times \Phi_*\eta) \bigr\} \\
	& = \mass(\xi; G) \cdot \mass(\eta; G)
\end{align*}
on homology level.
To show the inequality of the stable systoles,
recall that the cross product homomorphism
\begin{align*}
	H_p(X,A;\mathbb{R}) \otimes H_q(Y,B;\mathbb{R}) \to
	H_{p+q}\bigl( (X,A) \times (Y,B) ; \mathbb{R} \bigr)
\end{align*}
is a monomorphism.
Therefore we can estimate the stable $q$--systole as
\begin{align*}
	\stsys_{p+q}\bigl( (X,A) \times (Y,B) \bigr)
	& \le \inf \left\{\mass(\xi \times \eta) : 
	\begin{tabular}{l}
		$\xi \in H_p(X,A;\mathbb{Z})$, $\iota_*\xi \neq 0$, \\
		$\eta \in H_q(Y,B;\mathbb{Z})$, $\iota_*\eta \neq 0$
	\end{tabular}
	\right\} \\
	& \le \stsys_p(X,A) \cdot \stsys_q(Y,B) \,.
\end{align*}
where the second inequality is obtained by the result on homology level.
\end{proof}

\begin{lemma}
\label{lemma:StsysProjectionEquality}
Suppose $X$ and $Y$ are local Lipschitz neighborhood retracts.
If $Y$ is connected and the K\"{u}nneth formula gives an isomorphism of non-trivial vector spaces
\begin{align*}
	H_q(X;\mathbb{R}) \otimes H_0(Y;\mathbb{R}) \cong H_q\bigl( X \times Y ;\mathbb{R} \bigr) \neq \{0\} \,,
\end{align*}
then the stable $q$--systole satisfies
$$0 < \stsys_q\bigl( X \times Y \bigr) = \stsys_q(X) < \infty.$$
with respect to the product metric on $X \times Y$.
\end{lemma}

\begin{proof}
Let $\mathfrak{pr}_1\co X \times Y \to X$ be the first projection.
From the assumption,
for a nonzero homology class $\eta$ in $H_q ( X \times Y ;\mathbb{R} )$,
there exist $[S] \neq 0$ in $H^{\flat}_q(X;\mathbb{R})$ and $[T] \neq 0$ in $H^{\flat}_0(Y;\mathbb{R})$
whose cross product is the image of $\eta$ 
in $H^{\flat}_q\bigl(X \times Y;\mathbb{R}\bigr)$ with the same positive mass, ie,
$$\mass\bigl( [S] \times [T] \bigr) = \mass(\eta) > 0.$$
Note that the vector space of normal $0$--chains $\normalD_0(Y;\mathbb{R})$ is equal to the vector space of polyhedral $0$--chains $\polyhedralD_0(Y;\mathbb{R})$ which is generated by $\{\mathfrak{d}_y : y \in Y \}$ where $\mathfrak{d}$ is defined in the proof of \fullref{prop:stable0systole}.
For every points $y$ and $y'$ in $Y$,
$[\mathfrak{d}_y] = [\mathfrak{d}_{y'}]$ implies that
there is a nonzero real number $r$ such that
$[T] = r [\mathfrak{d}_y]$ with $\mass [T] = |r| \cdot \mathfrak{d}_y (1_Y^*) = |r|$.
Also, every $[S] \times [T]$ has representation of $[r \cdot S] \times [\mathfrak{d}_y]$,
therefore $\mathfrak{pr}_{1*}$ is an isomorphism with
$\mathfrak{pr}_{1*} \bigl( [S] \times [T] \bigr) = [r \cdot S]$.
Hence \fullref{lemma:StableSystoleMappedInequality} implies
$$\stsys_q(X \times Y) \ge \stsys_q(X) > 0$$
with the fact of $\mathfrak{pr}_{1}$ is a Lipschitzian map with $\Lip(\mathfrak{pr}_{1}) = 1$.
As a result, we obtain the equality by combining the result of \fullref{lemma:QuotientMassCrossEquality}.
\end{proof}

\section{Calculation by dimension and constructing metrics}
\label{section:Result}

At first, we will calculate the stable systolic category from the dimensional information of homology.
If the homology group is not so complex such as a real homology sphere,
we know the stable systolic category by only using dimensional information.
If a oriented manifold has a relatively simple cup-product structure such as $n$--fold producted space of spheres,
then the stable systolic category can be also calculated instantly.
Such methods to calculate the stable systolic category can be generalized as follows.

For a topological space $X$, let $\LPD(X)$ denote the {\it least positive dimension} of real cohomology vector spaces of $X$. 
So $\LPD(X) = l$ if and only if $\tilde{H}^i(X;\mathbb{R}) = \{0\}$ for $0 < i < l$ and $\tilde{H}^l(X;\mathbb{R}) \neq \{0\}$.
If $M$ is an $m$--manifold, then $\LPD(M)$ is smaller than $m$.

\begin{defn}
An $n$--dimensional CW space $X$ is said to {\it have maximal real cup length},
if there exist some real cohomology classes $\alpha_1, \cdots, \alpha_r$ with $\alpha_i \in \tilde{H}^{d_i}(X;\mathbb{R})$, a nonzero cup-product $\alpha_1 \cup \cdots \cup \alpha_r \in \tilde{H}^n(X;\mathbb{R})$ and $r := \lfloor n/\LPD(X) \rfloor$ where $\lfloor x \rfloor$ denotes the floor of a real number $x$.
\end{defn}

\begin{example}
Let $S$ be a manifold which is a real homology sphere.
Then $S$ has maximal real cup length, because of $\LPD(S) = \dim(S)$.
The $n$--fold direct product of $S$ also has maximal real cup length.
The direct product $S^2 \times S^3$ of spheres has maximal real cup length.
\end{example}

\begin{cor}
If an $m$--manifold $M$ has maximal real cup length,
then the stable systolic category of $M$ is equal to the real cup-length of $M$, ie,
$$\catstsys(M) = \cuplength_{\mathbb{R}}(M) = \lfloor m/\LPD(M) \rfloor .$$
\end{cor}

\begin{proof}
We need to verify that $\catstsys(M) \le \cuplength_{\mathbb{R}}(M)$.
Let $r := \lfloor m/\LPD(M) \rfloor$.
If $(d_1, \cdots, d_{k})$ is a partition of $m$ such that each stable $d_i$--systole is non-trivial, then $d_i \ge \LPD(M)$, so
 there is an inequality
$${k} \cdot \LPD(M) \le m = d_1+ \cdots + d_{k} < (r+1) \cdot \LPD(M)$$
which implies ${k} \le r = \cuplength_{\mathbb{R}}(M)$.
\end{proof}

In general, the direct product $M \times N$ of manifolds does not have maximal real cup length even if $M$ and $N$ have maximal real cup-length.
For example, the direct product of spheres $S^1 \times S^2$ does not have maximal real cup length.

\begin{lemma}\label{lemma:lowerboundP}
If manifolds $M_1^{m_1}, \cdots, M_n^{m_n}$ have maximal real cup length, then the stable systolic category of their $n$--fold direct product $M_1 \times \cdots \times M_n$ is greater than the sum of stable systolic categories for each $M_i$, ie,
$$\catstsys\left( M_1 \times \cdots \times M_n \right) \ge \catstsys(M_1) + \cdots + \catstsys(M_n).$$
\end{lemma}

\begin{proof}
Since $M_i$ has maximal real cup length,
there is nonzero cup product
$\bigcup_{j=1}^{r_i} \alpha_i^j$ in $H^{m_i}(M_i;\mathbb{R})$
where $r_i := \lfloor m_i/\LPD(M_i) \rfloor = \catstsys(M_i)$ for $1 \le i \le n$.

By the K\"{u}nneth formula, the $n$--fold cross product on the top dimensions induces an isomorphism
\begin{equation*}
	\bigotimes_{i=1}^n H^{m_i}(M_i;\mathbb{R}) 
	\cong H^{m}\left(\textstyle\prod_{i=1}^n M_i; \mathbb{R}\right) \quad \text{where} \quad m := \textstyle\sum\limits_{i=1}^n m_i .
\end{equation*}
This implies the existence of a nonzero cup product
\begin{align*}
	\prod_{i=1}^n \textstyle\left(\bigcup_{j=1}^{r_i} \alpha_i^j\right)
	& = \bigcup_{i=1}^n \textstyle \mathfrak{pr}_i^*\biggl(\,\bigcup\limits_{j=1}^{r_i} \alpha_i^j\biggr)
	 = \textstyle\bigcup\limits_{j=1}^{r_1} \mathfrak{pr}_1^*\alpha_1^j \cup \cdots \cup \bigcup\limits_{j=1}^{r_n} \mathfrak{pr}_n^*\alpha_n^j
\end{align*}
in the top-dimensional real cohomology vector space $H^{m}\left(\textstyle\prod_{i=1}^n M_i; \mathbb{R}\right)$, where $\mathfrak{pr}_i\co M_1 \times \cdots \times M_n \to M_i$ is the $i$--th projection.
This cup product implies that $\catstsys\bigl(\textstyle\prod_{i=1}^n M_i\bigr) \ge r_1 + \cdots + r_n$ by Gromov's Theorem.
\end{proof}

\begin{prop}\label{prop:LPDofProduct}
For manifolds $M$ and $N$, the least positive dimension of cohomology of $M \times N$ is
the minimum of $\LPD(M)$ and $\LPD(N)$.
\end{prop}

\begin{proof}
From the K\"{u}nneth formula, $H^i(M \times N;\mathbb{R}) = \{0\}$ for $0 < i < \min\bigr(\LPD(M),\LPD(N)\bigl)$.
If $l := \min\bigr(\LPD(M),\LPD(N)\bigl) = \LPD(M)$, then 
$H^l(M;\mathbb{R})$ is nonzero and
the cross product homomorphism $H^l(M;\mathbb{R}) \otimes H^0(N;\mathbb{R}) \to H^l(M \times N;\mathbb{R})$ is a monomorphism.
Therefore $H^l(M \times N;\mathbb{R})$ is nonzero.
The case of $\LPD(M) > \LPD(N)$ is shown by using the same arguments.
\end{proof}

For integers $i$ and $j \neq 0$,
let $\MOD(i,j)$ denotes the remainder from the division of $i$ by $j$.

\begin{cor}\label{cor:productedLSSCP}
Suppose manifolds $M^m$ and $N^n$ have maximal real cup length, and an integer 
$l := \LPD(M \times N)$.
If $M$ and $N$ satisfy the conditions
\begin{align*}
	& \lfloor m / \LPD(M) \rfloor = \lfloor m / l \rfloor, \quad\quad
	\lfloor n / \LPD(N) \rfloor = \lfloor n / l \rfloor \\
\tag*{\text{\sl and}}
	& \MOD(m,l) + \MOD(n,l) < l ,
\end{align*}
then $M \times N$ has maximal real cup length. Therefore,
$$\catstsys(M \times N) = \catstsys(M) + \catstsys(N) \,.$$
\end{cor}

\begin{proof}
Let integers $r := \lfloor m / l \rfloor$ and $s := \lfloor n / l \rfloor$.

\fullref{prop:LPDofProduct} implies that
$l = \min\bigl(\LPD(M),\LPD(N)\bigr) = \LPD(M \times N)$.
So we can formulate
$\lfloor (m+n) / \LPD(M \times N) \rfloor = r + s + \lfloor \MOD(m,l) + \MOD(n,l) \rfloor$. 
By the assumption, $\lfloor \MOD(m,\LPD(M)) + \MOD(n,\LPD(N)) \rfloor$ is zero, so we have
$$\lfloor (m+n) / \LPD(M \times N) \rfloor = r + s .$$
Thus it is sufficient to show that there is a nonzero cup product with the length of $r+s$.

Since $M$ and $N$ have maximal real cup length, there are cohomology classes $\alpha_1, \cdots, \alpha_r$ and $\beta_1, \cdots, \beta_s$ with their cup products are nonzero cohomology classes $\bigcup_{i=1}^r \alpha_i$ in $H^m(M;\mathbb{R})$ and $\bigcup_{i=1}^s \beta_i$ in $H^n(M;\mathbb{R})$.
From the proof of \fullref{lemma:lowerboundP},
there is a nonzero cup product
$\bigcup_{i=1}^r \mathfrak{pr}_1^*\alpha_i \cup \bigcup_{i=1}^s \mathfrak{pr}_2^*\beta_i$ in the top dimensional cohomology vector space $H^{m+n}(M \times N;\mathbb{R})$.
\end{proof}

Without the condition of the product $M \times N$ has maximal real cup length, we can generalize this corollary as follow.

\begin{thm}
\label{thm:LSSCPstsyscat}
Let manifolds $M^m$ and $N^n$ have maximal real cup length.
If 
$$\MOD\bigl(m,\LPD(M)\bigr) + \MOD\bigl(n,\LPD(N)\bigr) < \max\bigl(\LPD(M),\LPD(N)\bigr) \,,$$
then the stable systolic category of their product $M \times N$ is the sum of each stable systolic category, ie,
$$\catstsys(M \times N) = \catstsys(M) + \catstsys(N) \,.$$
\end{thm}

\begin{proof}
Since $M$ and $N$ have maximal real cup length,
$$r := \lfloor m/\LPD(M) \rfloor = \catstsys(M) \quad \text{ and } \quad s := \lfloor n/\LPD(N) \rfloor = \catstsys(N) \,.$$
In the case of $\LPD(M) = \LPD(N)$ is \fullref{cor:productedLSSCP}.
So we will assume $\LPD(M) < \LPD(N)$.

From \fullref{lemma:lowerboundP}, 
$\catstsys(M \times N) \ge \catstsys(M) + \catstsys(N) = r+s$ .
Therefore, it is sufficient to show that any partition of $m{+}n$ whose size is greater than $r{+}s$, is not a stable systolic categorical partition.

Suppose $m+n = d_1 + \cdots + d_k$ is a stable systolic categorical partition for $M \times N$ with some integer $1 \le r' \le k$ and the condition $0 < \LPD(M) \le d_1 \le \cdots \le d_{r'} < \LPD(N)$.
For an arbitrary $t \ge 1$, let $\mathcal{G}_{t} := t^2 \mathcal{G}_M + \mathcal{G}_N$ be a Riemannian metric on $M \times N$.
Then \fullref{lemma:MassRescaleEquality} and \fullref{lemma:StsysProjectionEquality}
imply that the stable systoles for the partition $d_1 + \cdots + d_k$ satisfies
\begin{align*}
	\prod_{i=1}^k\stsys_{d_i}(M \times N, \mathcal{G}_{t}) 
	& \ge \prod_{i=1}^{r'}\stsys_{d_i}(M, t^2\mathcal{G}_{M}) \prod_{j=r'+1}^{k}\stsys_{d_j}(M \times N, \mathcal{G}_{t}) \\
	& = t^{d_1 + \cdots + d_{r'}}\prod_{i=1}^{r'}\stsys_{d_i}(M, \mathcal{G}_{M}) \prod_{j=r'+1}^{k}\stsys_{d_j}(M \times N, \mathcal{G}_{t})
\end{align*}
Since we assume that $t \ge 1$, we can obtain the inequality
$\stsys_{d_j}(M \times N, \mathcal{G}_{t}) \ge \stsys_{d_j}(M \times N, \mathcal{G}_{1})$
for each $r'+1 \le j \le k$.
On the other hands,
the mass of integral fundamental class $[M \times N]$ is characterized by 
\fullref{lemma:MassRescaleEquality} and \fullref{lemma:QuotientMassCrossEquality} as
\begin{align*}
	\mass\bigl( [M \times N],\mathcal{G}_t \bigr)
	& \le \mass\bigl( [M],t^2\mathcal{G}_M \bigr) \cdot \mass\bigl( [N], \mathcal{G}_N \bigr) \\
	& = t^m \cdot \mass\bigl( [M], \mathcal{G}_M \bigr) \cdot \mass\bigl( [N], \mathcal{G}_N \bigr) .
\end{align*}
Here if we assume that $d_1 + \cdots + d_{r'} > m$, then we have
\begin{align*}
	\frac{\prod\limits_{i=1}^k\stsys_{d_i}(M {\times} N, \mathcal{G}_{t})}{\mass\bigl( [M {\times} N],\mathcal{G}_t \bigr)}
	& \ge t^{m - (d_1 + \cdots + d_{r'})} \cdot \frac{\prod\limits_{i=1}^{r'}\stsys_{d_i}(M, \mathcal{G}_{M}) \cdot \prod\limits_{j=r'+1}^{k}\stsys_{d_j}(M {\times} N, \mathcal{G}_{1})}{\mass\bigl( [M], \mathcal{G}_M \bigr) \cdot \mass\bigl( [N], \mathcal{G}_N \bigr)}
\end{align*}
where the right-hand side of the inequality diverges as $t \to \infty$. 
This contradicts to that $(d_1, \cdots, d_k)$ is a stable systolic categorical partition.
Hence we obtain $d_1 + \cdots + d_{r'} \le m$ and $d_{r'+1} + \cdots + d_{k} \ge n$.
This condition for $m$ implies
$$r' \le \lfloor (d_1 + \cdots + d_{r'}) / \LPD(M) \rfloor \le \lfloor m / \LPD(M) \rfloor \le r \,.$$
Let $s' := k - r'$.
From the assumption, $\LPD(M) / \LPD(N) < 1$ and
$$\MOD(m,\LPD(M)) + \MOD(n,\LPD(N)) < \LPD(N),$$
so we can calculate as
\begin{align*}
	& k = r' + s' \le r + s
\end{align*}
which implies $\catstsys(M \times N) \le \catstsys(M) + \catstsys(N)$ .
\end{proof}

\begin{cor}
Suppose manifolds $M_0 \times M_1 \times \cdots \times M_{k}$ and $M_{k+1} \times \cdots \times M_n \times M_{n+1}$
have maximal real cup length with 
\begin{align*}
	& \LPD(M_0) = \LPD(M_1) = \cdots = \LPD(M_k) \\
\tag*{\text{\sl and}}
	& \LPD(M_{k+1}) = \cdots = \LPD(M_n) = \LPD(M_{n+1}) \,.
\end{align*}
Let $r_i := \lfloor \dim(M_i)/\LPD(M_i) \rfloor$ for $0 \le i \le n+1$ .
If $M_0, \cdots, M_{n+1}$ satisfy conditions
$\dim(M_i) = \LPD(M_i) \cdot r_i$ for $1 \le i \le n$ and
\begin{eqnarray*}
	\dim(M_0) - \LPD(M_0)\cdot r_0 + \dim(M_{n+1}) - \LPD(M_{n+1})\cdot r_{n+1} \\
	< \max\bigr(\LPD(M_0),\LPD(M_{n+1})\bigl) ,
\end{eqnarray*}
\begin{gather*}
\tag*{\text{\sl then:}}
	\catstsys\left( {\textstyle\prod\limits_{i=0}^{n+1}M_i} \right) = \sum_{i=0}^{n+1} \catstsys(M_i) = \sum_{i=0}^{n+1} r_i 
\end{gather*}
\end{cor}

For the product $S^1 \times S^2$ of spheres, \fullref{thm:LSSCPstsyscat} can not applied.
So we must approach from the other viewpoints to obtain the stable systolic category.

\begin{thm}
\label{thm:CatstsysProductSphere}
If manifolds $S_1^{m_1}, \cdots, S_n^{m_n}$ are real homology spheres, then the stable systolic category of their $n$--fold direct product is the number of spheres.
\end{thm}

\begin{proof}
Since every real homology spheres have maximal real cup length, \fullref{lemma:lowerboundP}
gives us a lower bound 
$\catstsys(S_1 \times \cdots \times S_n) \ge n$.

Suppose $m_i \le m_{i+1}$ for each $1 \le i \le n$.
Then a partition $(m_1, \cdots, m_n)$ of $\sum_i m_i$ can be rewritten as $(r_1, \cdots, r_1, r_2, \cdots, r_{l-1}, r_l, \cdots, r_l)$ 
where $r_i$ is a range.
This corresponding to rewrite
\begin{align*}
	S_1^{m_1} \times \cdots \times S_n^{m_n} 
	& = \left(S_1^{r_1} \times \cdots \times S_{s_1}^{r_1}\right) \times \left(S_{s_1+1}^{r_2} \times \cdots \times S_{s_1+s_2}^{r_2}\right) \times \cdots \\
	& \qquad \times \left(S_{s_1+\cdots+s_{l-1}+1}^{r_l} \times \cdots \times S_{s_1+\cdots+s_{l-1}+s_l}^{r_l}\right)
\end{align*}
where $r_i := m_{s_1+\cdots+s_{i-1}+1} = \cdots = m_{s_1+\cdots+s_{i-1}+s_i}$ with $r_i < r_{i+1}$
and $s_i > 0$ is the duplicated number of $r_i$, so that $s_1 + \cdots + s_l = n$.
For simplicity, let
$$X_p := S_1 \times \cdots \times S_{s_1 + \cdots + s_p}
\quad \text{and} \quad
Y_p := S_{s_1 + \cdots + s_p} \times \cdots \times S_{n}$$
for $1 \le p \le n$.
Then $S_1 \times \cdots \times S_n = X_p \times Y_p$
and we can observe that
$\mathcal{G}_{p,t} := t^2 \mathcal{G}_{X_p} + \mathcal{G}_{Y_p}$
is a Riemannian metric on $X_p \times Y_p$ for $t > 0$ when $\mathcal{G}_{X_p} + \mathcal{G}_{Y_p}$ is a Riemannian metric on $X_p \times Y_p$.
Now we can apply \fullref{lemma:StsysProjectionEquality} and \fullref{lemma:MassRescaleEquality}, so
$$\stsys_q(X_p \times Y_p, \mathcal{G}_{p,t}) = \stsys_q(X_p, t^2 \mathcal{G}_{X_p}) = t^q \cdot \stsys_q(X_p, \mathcal{G}_{X_p})$$
for the non-trivial stable systoles in the dimension of $1 \le q \le s_1 + \cdots + s_p$.

Here we suppose $(d_1, \cdots, d_k)$ with $d_i \le d_{i+1}$, is the longest stable systolic categorical partition for $S_1 \times \cdots \times S_n$.
Then we can rewrite $(d_1, \cdots, d_k)$ by the ranges $\{r_1, \cdots, r_l\}$
and the duplicated number $s_i' \ge 0$ of $r_i$.
We will show this partition is not longer than $n$ by induction on $p$ for $1 \le p \le l$ and contradiction.
Assume that $s_i' = s_i$ for $1 \le i \le p-1$.
If $s_p' > s_p$, then
using a similar argument in the proof of \fullref{thm:LSSCPstsyscat}, we can observe that the right-hand side of the inequality
\begin{align*}
	\frac{\prod\limits_{i=1}^{k}\stsys_{d_i}(X_p {\times} Y_p, \mathcal{G}_{p,t})}{\mass\bigl([X_p {\times} Y_p], \mathcal{G}_{p,t}\bigr)}
	& \ge t^{w} \cdot \frac{\prod\limits_{i=1}^{p}\stsys_{r_i}(X_p, \mathcal{G}_{X_p})^{s_i'} \prod\limits_{i=p+1}^{l}\stsys_{r_i}(X_p {\times} Y_p, \mathcal{G}_{p,1})^{s_i'}}{\mass\bigl([X_p], \mathcal{G}_{X_p}\bigr) \cdot \mass\bigl([Y_p], \mathcal{G}_{Y_p}\bigr)}
\end{align*}
diverge as $t \to \infty$ where 
$w := r_1(s'_1 - s_1) + \cdots + r_i(s'_p - s_p) > 0$.
This contradicts to that the partition $(d_1, \cdots, d_k)$ is stable systolic categorical,
and hence we obtain $s_p' \le s_p$.
However we must choose $s_p' = s_p$ to make the longest partition.
As a result, the size of the longest stable systolic categorical partition can not exceed $n = s_1 + \cdots + s_l$.
\end{proof}

\section{Invariance under the rational equivalences}
\label{section:RHEinvariant}

Suppose $M$ and $N$ are $n$--manifolds.
Let $K$ and $L$ be a triangulation of $M$ and $N$ respectively.
In this section, $K$ and $L$ are subdivided if necessary, but we will use the same symbol.
For a continuous map $f\co M \to N$, 
there is a non-degenerate simplicial approximation $g\co K \to L$ of $f$.
For an open $n$--simplex $e$ in $L$, consider a map $h\co K \overset{g}{\to} L \to L / (L \setminus e)$.
We will call $\deg(h)$ the {\it degree} of $g$ at $e$ which is denoted by $\deg_e(g)$.
Let 
$$D(g) := \sup\bigl\{ \vert \deg_e(g) \vert : \text{open } n \text{-simplex } e \text{ in } L \bigr\}.$$
Here $D(g)$ is finite, because of we can assume that $K$ and $L$ are finite simplicial complexes.

For an arbitrary Riemannian metric $\mathcal{G}_N$ on $N$ and $\varepsilon > 0$, 
there is a piecewise linear metric $\mathcal{G}_L = \mathcal{G}_L(\varepsilon)$ on $L$ satisfying
$$\bigl|\, \stsys_q(L,\mathcal{G}_L) - \stsys_q(N,\mathcal{G}_N) \,\bigr| \le \varepsilon$$
for every non-trivial stable $q$--systoles (compare Federer \cite[4.1.22]{Fed69}) and the realization of $L$ with $\mathcal{G}_L$ is a PL section of the normal bundle over $N$ with $\mathcal{G}_N$ in some finite dimensional Euclidean space.
Such metric can be obtained by subdividing $K$ and $L$, and translating vertices in $L$ along the fiber of the normal bundle to do not degenerate any simplex.
For $0 < \varepsilon' < \varepsilon$,
a suitable metric $\mathcal{G}_L(\varepsilon')$ also can be acquired by the same way.
Hence we can assume that $D(g)$ is not changed by $\varepsilon$ and $\mathcal{G}_L$.
As $\varepsilon$ approaches to $0$, each $L$, $\mathcal{G}_L$ and $g^*\mathcal{G}_L$ converges to $N$, $\mathcal{G}_N$ and a piecewise Riemannian metric on $M$ respectively.
Under this circumstance, we obtain following lemma.

\begin{lemma}
\label{lemma:EquivalentStableSystole}
Suppose $q$--th real homology vector space of $K$ and $L$ are non-trivial.
If $g\co K \to L$ induces a monomorphism $g_*$ between the $q$--th real homology vector spaces, then
$$\stsys_q(L,\mathcal{G}_L) \le \stsys_q(K,g^*\mathcal{G}_L) \le D(g) \cdot \stsys_q(L,\mathcal{G}_L) < \infty$$
for every piecewise linear metric $\mathcal{G}_L$ on $L$.
\end{lemma}

\begin{proof}
With the pullback PL metric $g^*\mathcal{G}_L$ on $K$,
$g$ is a distance decreasing map.
Combining this with \fullref{lemma:StableSystoleMappedInequality},
$$\stsys_q(L,\mathcal{G}_L) \le \Lip(g)^q \cdot \stsys_q(K,g^*\mathcal{G}_L) \le \stsys_q(K,g^*\mathcal{G}_L).$$

On the other hands, 
the inverse image of an arbitrary $q$--simplex of $L$ is $D(g)$ of $q$--simplices as at most,
since $g$ is a non-degenerate simplicial map and
every $q$--simplex is contained in the boundary of some $n$--simplex for $q < n$.
Also each simplex in the inverse image has same mass of the preimage, since the restriction of $g$ on each simplex is isometry. 
This implies that the mass of a $q$--chain $c$ of $K$ is not greater than $D(g)$ times of the mass of the image $g_{\flat}(c)$ which is not trivial.
Therefore we can verify that
$$\stsys_q(K,g^*\mathcal{G}_L) \le D(g) \cdot \stsys_q(L,\mathcal{G}_L)$$
for an arbitrary PL metric $\mathcal{G}_L$.
\end{proof}

\begin{rmk}
If $K$ is not a triangulation of a manifold,
we can not sure that every $q$--simplex of $K$ is contained in the boundary of some $n$--simplex for $q < n$.
For example, a triangulation of the one-point union $S^1 \vee S^2$ has some $1$--simplex in $S^1$ which is not contained in the boundary of any $2$--simplex.
\end{rmk}

Since the stable systolic category is a homotopy invariant,
here we obtain following proposition using similar techniques of Katz and Rudyak \cite{KatRud08}.

\begin{prop}
Let $M$ and $N$ are $n$--manifolds.
If there exists a smooth map $f\co M \to N$ which induces a monomorphism on every real homology vector space, then
$\catstsys(M) \le \catstsys(N)$.
\end{prop}

\begin{proof}
We apply \fullref{lemma:EquivalentStableSystole},
\begin{align*}
	& \stsys_{q}(N,\mathcal{G}_N) 
	\le \stsys_{q}(L,\mathcal{G}_L) + \varepsilon
	\le \stsys_{q}(K,g^*\mathcal{G}_L) + \varepsilon \\
\tag*{\text{and}}
	& \stsys_{q}(N,\mathcal{G}_N) + \varepsilon
	\ge \stsys_{q}(L,\mathcal{G}_L)
	\ge 1/D(g) \cdot \stsys_{q}(K,g^*\mathcal{G}_L)
\end{align*}
where
$L$ converges to $N$ in some Euclidean space and 
$g^*\mathcal{G}_L$ converges to a piecewise Riemannian metric $\mathcal{G}_M$ on $M$ as $\varepsilon$ approaches to $0$.
Suppose there exists a stable systolic categorical partition $(d_1, \cdots, d_k)$ for $M$.
Then there exist $C > 0$ and $\delta = \delta(\varepsilon) > 0$ such that $\delta$ converges to $0$ as $\varepsilon$ approaches to $0$ and
$$\prod_{i=1}^k \stsys_{d_i}(K,g^*\mathcal{G}_L)
\le C \cdot \mass([K],g^*\mathcal{G}_L) + \delta,$$
because of each metric $g^*\mathcal{G}_L$ can be approximated by some Riemannian metrics on $M$.
We can assume that $\varepsilon \le \stsys_{d_i}(N,\mathcal{G}_N)$ for all $i$, so
\begin{align*}
	\prod_{i=1}^k \stsys_{d_i}(L,\mathcal{G}_L) 
	& \le 2^k \cdot \prod_{i=1}^k \stsys_{d_i}(K,g^*\mathcal{G}_L)  \\
	& \le 2^k C \cdot \mass([K],g^*\mathcal{G}_L) + 2^k \delta \\
	& \le 2^k C A(g) \cdot \mass([L],\mathcal{G}_L) + 2^k(C A(g) \cdot \varepsilon + \delta).
\end{align*}
This implies the partition $(d_1, \cdots, d_k)$ is also stable systolic categorical for $N$.
Therefore we obtain the result $\catstsys(M) \le \catstsys(N)$.
\end{proof}


Let $X$ and $Y$ are simply connected spaces.
A continuous map $f\co X \to Y$ is called a {\it rational equivalence}, 
if the induced map $f^*\co H^*(Y; \mathbb{Q}) \to H^*(X; \mathbb{Q})$ is an isomorphism.

\begin{cor}
\label{cor:REinvariance}
The stable systolic category of a $0$--universal manifold is invariant under the rational equivalences.
\end{cor}

\begin{proof}
For a $0$--universal manifold $M$ and a rational equivalence to a space $X$, there exists a rational equivalence from $X$ to $M$.
\end{proof}

\bibliographystyle{gtart}
\bibliography{reference.paper}

\end{document}